\documentclass[12pt]{amsart}
\usepackage{amsmath,graphicx}
\usepackage[all]{xy}

\newtheorem{theorem}{Theorem}         

\newtheorem{remark}{Remark}

\providecommand{\area}{\mathop{\rm area}\nolimits}
\providecommand{\Real}{\mathop{\rm Re}\nolimits}%
\providecommand{\supp}{\mathop{\rm supp}\nolimits}
\providecommand{\R}{{\mathbb{R}}}
\providecommand{\FF}{{\mathcal{F}}}
\providecommand{\BB}{\mathcal{B}}
\providecommand{\Z}{{\mathbb{Z}}}
\providecommand{\RR}{{\mathcal{R}}}
\providecommand{\N}{\mathbb{N}} \providecommand{\eps}{\varepsilon}
\providecommand{\SSS}{\mathcal{S}}

\numberwithin{equation}{section}

\begin{document}

\title[]{The correlations of Farey fractions}

\author[]{Florin P. Boca and Alexandru Zaharescu}

\thanks{Research partially supported by the Ceres
Programme-competition C4}

\address{\sl Department of Mathematics, University of Illinois,
1409 W. Green Str., Urbana IL 61801, USA}

\address{\sl Institute of Mathematics of the Romanian Academy,
P.O. Box 1-764, RO-014700, Romania}

\address{E-mail: fboca@math.uiuc.edu; zaharesc@math.uiuc.edu}

\date{January 18, 2005}

\begin{abstract}
We prove that all correlations of the sequence of Farey fractions
exist and provide formulas for the correlation measures.
\end{abstract}

\maketitle

\section{Introduction}

By the classical contributions of Franel and Landau (\cite{Fra},
\cite{Lan}, see also \cite{Edw}), the existence of zero-free
regions $1-\delta_0<\Real s<1$ for the Riemann zeta function - and
in particular the Riemann hypothesis - are known to be equivalent
to quantitative statements about the uniform distribution of Farey
fractions. Although these problems remain widely open, the study
of the distribution of Farey fractions is of independent interest.
Ideas and techniques from this area, especially Weil-Sali\' e type
estimates on Kloosterman sums, turned out to be useful in proving
sharp asymptotic formulas in problems from various areas of
mathematics (\cite{BCZ}, \cite{BGZ}, \cite{BZ}, \cite{HZ},
\cite{Sar}).

The statistics of the spacings between the elements of a (scaled)
sequence of real numbers can be expressed in the convergence of
certain measures, called consecutive spacing measures and
respectively correlation measures. The first ones reflect the
distribution of consecutive tuples of elements. Correlations
measure the distribution between all tuples of elements and do not
require the ordering of the numbers. There are very few sequences
of interest for which one could establish the existence of
correlation measures and many of them are conditional, as in the
important case of the zeros of the Riemann zeta function, or more
general $L$-functions (\cite{Mo}, \cite{Hej}, \cite{RS},
\cite{KS}). In this note we study the correlations of the sequence
$\FF_Q$ of Farey fractions of order $Q$ in $[0,1]$ as
$Q\rightarrow \infty$.

Let $\nu \geq 1$ be an integer and let $F$ be a finite set of $N$
elements in $[0,1]$. The \emph{$\nu$-level correlation measure}
$\RR_F^{(\nu)} (\BB)$ of a box $\BB \subset {\mathbb{R}}^{\nu -1}$
is defined as
\begin{equation*}
\frac{1}{N}\, \# \Big\{ (x_1,\dots,x_\nu) \in F^\nu  : x_i \
\mbox{\rm distinct},\, (x_1-x_2,x_2-x_3,\dots,x_{\nu -1}-x_\nu)\in
\frac{1}{N}\ \BB +\Z^{\nu -1}\Big\}.
\end{equation*}
When $\nu=2$, the \emph{pair correlation measure} of an interval
$I\subset \R$ is
\begin{equation*}
\RR^{(2)}_F (I)=\frac{1}{N}\, \# \Big\{ (x,y)\in F^2  : x\neq y\
\mbox{\rm and}\ x-y\in \frac{1}{N} \ I+\Z\Big\} .
\end{equation*}
Suppose that $(F_n)_n$ is an increasing sequence of finite subsets
of $[0,1]$ and that
\begin{equation*}
\RR^{(\nu)}(\BB)=\lim\limits_n \RR^{(\nu)}_{F_n}(\BB)
\end{equation*}
exists for every box $\BB \subset \RR^{\nu -1}$. Then
$\RR^{(\nu)}$ is called the \emph{$\nu$-level correlation measure
of $(F_n)_n$}. The measure $\RR^{(2)}$ is called the \emph{pair
correlation measure of $(F_n)_n$.} If
\begin{equation*}
\RR^{(\nu)}(\BB)=\int\limits_\BB g_\nu (x_1,\dots,x_{\nu -1}) dx_1
\dots dx_{\nu -1},
\end{equation*}
then $g_\nu$ is called the \emph{$\nu$-level correlation function
of $(F_n)_n$}. We denote
\begin{equation*}
\RR^{(\nu)}_F (\lambda_1,\dots,\lambda_{\nu-1})=2^{-\nu+1}
\RR^{(\nu)}_F \left( \ \prod\limits_{j=1}^{\nu -1}
[-\lambda_j,\lambda_j]\right) .
\end{equation*}

Given $\Lambda >0$ and two $(\nu-1)$-tuples
$A=(A_1,\dots,A_{\nu-1}),\ B =(B_1,\dots,B_{\nu -1})\in \N^{\nu
-1}$, we consider the map defined by
\begin{equation}\label{1.1}
T_{A,B}(x,y)=\frac{3}{\pi^2} \left( \frac{B_1 }{y(A_1y-B_1
x)},\dots ,\frac{B_{\nu-1}}{y(A_{\nu-1}y-B_{\nu-1}x)} \right)
\end{equation}
and the set
\begin{equation}\label{1.2}
\Omega_{A,B,\Lambda} =\left\{ (x,y) : 0\leq x\leq y\leq 1, \ y\geq
\frac{3}{\pi^2 \Lambda},\ yA-xB\in (0,1]^{\nu-1} \right\}.
\end{equation}
We also consider the linear transformation on $\R^{\nu-1}$ defined
by
\begin{equation}\label{1.3}
T(x_1,x_2,\dots,x_{\nu-1})=
(x_1-x_2,x_2-x_3,\dots,x_{\nu-2}-x_{\nu-1},x_{\nu-1}),
\end{equation}
whose inverse is
\begin{equation*}
T^{-1}(y_1,\dots,y_{\nu-1})=(y_1+\cdots+y_{\nu-1},
y_2+\cdots+y_{\nu-1},\dots,y_{\nu-2}+y_{\nu-1},y_{\nu-1}),
\end{equation*}
and set
\begin{equation*}
\Phi_{A,B}=T\circ T_{A,B}.
\end{equation*}

The main result of this paper shows the existence of correlation
measures for the sequence $(\FF_Q)_Q$ of Farey fractions. This
complements previous results (\cite{Hall}, \cite{ABCZ}) concerning
the existence and computation of consecutive spacing measures
between Farey fractions and gives a complete description of their
spacing statistics.

\begin{theorem}\label{T1}
All $\nu$-level correlation measures of the sequence $(\FF_Q)_Q$
exist. Moreover, for any box $\BB \subset (0,\Lambda)^{\nu-1}$,
\begin{equation}\label{1.4}
\RR^{(\nu)}(\BB)=2\sum\limits_{\substack{A,B\in \N^{\nu-1} \\
(A_j,B_j)=1\\ j=1,\dots,\nu-1}} \area \big( \Omega_{A,B,\Lambda}
\cap \Phi_{A,B}^{-1} (\BB)\big).
\end{equation}
\end{theorem}

\begin{remark}
(i) For such a compact box $\BB$, $T^{-1}\BB$ is a compact set
that does not contain the origin. Since $0<y(A_jy-B_jx)\leq 1$ for
all $(x,y)\in \Omega_{A,B,\Lambda}$, all denominators in
\eqref{1.1} are greater or equal than $1$; thus the sum in
\eqref{1.4} has finitely many nonzero terms.

(ii) Since $(A_j,B_j)=1$, $j=1,\dots,\nu-1$, the maps $T_{A,B}$
and $\Phi_{A,B}$ are one-to-one.

(iii) For $\nu\geq 4$, the support of $\RR^{(\nu)}$ has Lebesgue
zero measure as a subset of a countable union of surfaces in
$\R^{\nu-1}$.
\end{remark}

When $\nu=2$ one gets a more explicit expression of the pair
correlation.

\begin{theorem}\label{T2}
The pair correlation function of $(\FF_Q)_Q$ is given by
\begin{equation}\label{1.5}
g_2 (\lambda)=\frac{6}{\pi^2 \lambda^2} \sum\limits_{1\leq k<
\frac{\pi^2\lambda}{3}} \varphi (k) \log \frac{\pi^2 \lambda}{3k}
.
\end{equation}
Moreover, as $\lambda \rightarrow \infty$,
\begin{equation}\label{1.6}
g_2 (\lambda)=1+O(\lambda^{-1}).
\end{equation}
\end{theorem}


The shape of the graph of $g_2$ shows some similarity with the
pair correlation function $g_{\mathrm{GUE}}(\lambda)=1-\sin^2 \pi
\lambda /(\pi^2 \lambda^2)$. Both distributions show repulsion
between the elements of the sequence. The strongest one occurs for
$g_2$ and is reflected by the vanishing of $g_2$ on the whole
interval $[0,3/\pi^2]$.

For a random (Poisson, uncorrelated) sequence of points, the pair
correlation function $g_{\mathrm{Po}}$ is identically equal to $1$
on $(0,\infty)$. As $\lambda \rightarrow \infty$, both pair
correlation functions $g_2$ and $g_{\mathrm{GUE}}$ approach $1$.

\begin{figure}[ht]
\begin{center}
\includegraphics*[scale=0.8, bb=0 0 300 200]{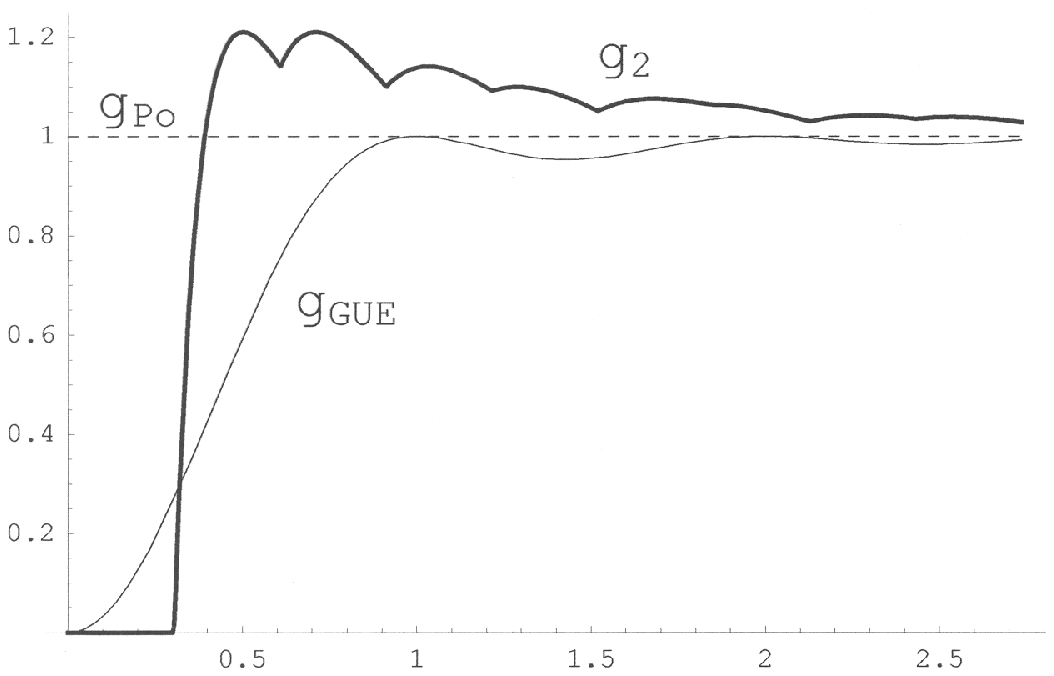}
\end{center}
\caption{\sl The graphs of $g_2$, $g_{\mathrm{GUE}}$, and
$g_{\mathrm{Po}}$.} \label{Figure1}
\end{figure}

\section{An expression of the smooth correlation sums}

It is well-known that the cardinality of $\FF_Q$ is
\begin{equation}\label{2.1}
N=N_Q=\sum\limits_{k=1}^Q \varphi (k)=\frac{3Q^2}{\pi^2}+O(Q\log
Q).
\end{equation}
An application of the Korobov-Vinogradov exponential sum estimates
gives the better estimate (see \cite{Wal})
\begin{equation}\label{2.2}
N_Q=\frac{3Q^2}{\pi^2}+O\big( Q(\log Q)^{2/3} (\log \log
Q)^{4/3}\big).
\end{equation}

If $\gamma=\frac{a}{q}$ and
$\gamma^\prime=\frac{a^\prime}{q^\prime}$ are two distinct
elements in $\FF_Q$, then $\vert \gamma^\prime -\gamma \vert \geq
\frac{1}{qq^\prime} \geq \frac{1}{Q^2}$. It now follows from
\eqref{2.1} that if $\min \vert \lambda_j \vert <\frac{3}{\pi^2}$,
then $\RR^{(\nu)} (\lambda_1,\dots,\lambda_{\nu -1})=0$. We also
have for any box $\BB \subset \RR^{\nu -1}$
\begin{equation*}
\RR^{(\nu)}(\BB)=\RR^{(\nu)} \big(\BB \setminus \{
(\lambda_1,\dots,\lambda_{\nu -1}) : \min \vert \lambda_j \vert <
3/\pi^2 \}\big).
\end{equation*}

We denote the M\" obius function by $\mu$ and the least common
multiple $\frac{d_1 d_2}{(d_1,d_2)}$ of two positive integers
$d_1$ and $d_2$ by $[d_1,d_2]$. We also set
\begin{equation*}
e(t)=\exp (2\pi it),\qquad t\in \R,
\end{equation*}
and
\begin{equation*}
x\cdot y=x_1y_1+\cdots+x_n y_n ,\qquad x=(x_1,\dots,x_n),\
y=(y_1,\dots,y_n)\in \R^n.
\end{equation*}
Exponential sums of the form
\begin{equation*}
\sum\limits_{\gamma \in \FF_Q} e(r\gamma),\qquad r\in \Z,
\end{equation*}
will play an important role in our computations. They are
connected with the function
\begin{equation*}
M(x)=\sum\limits_{n\leq x} \mu (n)
\end{equation*}
by the formula (cf. formula (1) at page 264 in \cite{Edw})
\begin{equation}\label{2.3}
\sum\limits_{\gamma \in \FF_Q} e(r\gamma) =
\sum\limits_{\substack{d\geq 1\\ d\vert r}} dM\bigg(
\frac{Q}{d}\bigg), \qquad r,Q\in \Z,\ Q\geq 1.
\end{equation}

We fix an integer $\nu \geq 2$, a smooth real-valued function $H$
on $\R^{\nu-1}$ such that $\supp (H)\subset (0,\Lambda)^{\nu -1}$,
and $0<\Lambda^\prime <\Lambda$ such that $\supp (H)\subset
(0,\Lambda^\prime)^{\nu-1}$. The Fourier transform of $H$ is
defined by
\begin{equation*}
\widehat{H}(x)=\int\limits_{\R^{\nu-1}} H(y) e(-x\cdot y)
dy,\qquad x\in \R^{\nu-1}.
\end{equation*}

We consider the $\Z^{\nu-1}$-periodic function $f$ given by
\begin{equation*}
f(y)=f_Q(y)=\sum\limits_{r\in \Z^{\nu -1}} H\big(
N(y+r)\big),\qquad y\in \R^{\nu -1} ,
\end{equation*}
and the smooth $\nu$-level correlation sum defined by
\begin{equation}\label{2.4}
\SSS_\nu =R^{(\nu)}(Q,H)=\frac{1}{N}
\sum\limits_{\gamma_1,\dots,\gamma_\nu \in \FF_Q \,\mathrm
{distinct}}\hspace{-5pt}
f(\gamma_1-\gamma_2,\gamma_2-\gamma_3,\dots,
\gamma_{\nu-1}-\gamma_\nu ).
\end{equation}

The Fourier coefficients in the Fourier series
\begin{equation*}
f(y)=\sum\limits_{r\in \Z^{\nu-1}} c_r  e(r\cdot y)
\end{equation*}
of $f$ are given by
\begin{equation}\label{2.5}
\begin{split}
c_r & =\int\limits_{[0,1)^{\nu-1}} f(y)e(-r\cdot y) dy
=\int\limits_{[0,1)^{\nu-1}} e(-r\cdot y) \sum\limits_{n\in
\Z^{\nu-1}} H\big( N(y+n)\big)dy \\ & =\sum\limits_{n\in
\Z^{\nu-1}} \int\limits_{[0,1)^{\nu-1}} e(-r\cdot y)  H\big(
N(y+n)\big) dy
\\ & =\sum\limits_{n\in \Z^{\nu-1}} \int\limits_{n+[0,1)^{\nu-1}}
\hspace{-10pt} e\big( -r\cdot (u-n)\big) H(Nu) du \\
& =\int\limits_{\R^{\nu-1}} e(-r\cdot u) H(Nu) du
=\frac{1}{N^{\nu-1}} \int\limits_{\R^{\nu-1}} e\bigg(
-\frac{r\cdot y}{N}\bigg) H(y) dy \\
& =\frac{1}{N^{\nu-1}} \widehat{H} \bigg( \frac{1}{N} \ r\bigg).
\end{split}
\end{equation}
Since $H$ is supported on $(0,\infty)^{\nu-1}$ we can remove in
\eqref{2.4}, for $Q$ large enough that $N>\Lambda$, the condition
that $\gamma_1,\dots,\gamma_\nu$ are distinct, and gather
\begin{equation}\label{2.6}
\begin{split}
& \SSS_\nu
=\frac{1}{N}\sum\limits_{\substack{\gamma_1,\dots,\gamma_\nu\in
\FF_Q \\ r_1,\dots,r_{\nu-1}\in \Z}} \hspace{-6pt} H\big(
N(r_1+\gamma_1-\gamma_2,r_2+\gamma_2-\gamma_3,\dots,r_{\nu-1}+
\gamma_{\nu-1}-\gamma_\nu)\big) \\
& \quad =\frac{1}{N} \sum\limits_{\gamma_1,\dots,\gamma_\nu \in
\FF_Q} f(\gamma_1-\gamma_2,\gamma_2-\gamma_3,\dots,\gamma_{\nu-1}
-\gamma_\nu)\\
& \quad =\frac{1}{N}
\sum\limits_{\substack{\gamma_1,\dots,\gamma_\nu\in \FF_Q \\
r_1,\dots,r_{\nu-1}\in \Z}} \hspace{-6pt} c_r \, e\big( r\cdot
(\gamma_1-\gamma_2,\gamma_2-\gamma_3,\dots,\gamma_{\nu-1}-\gamma_\nu)\big)
\\ &  =
\frac{1}{N} \sum\limits_{\substack{\gamma_1,\dots,\gamma_\nu\in
\FF_Q \\ r_1,\dots,r_{\nu-1}\in \Z}} \hspace{-6pt} c_r \,
e(r_1\gamma_1) e\big( (r_2-r_1)\gamma_2\big) \cdots e\big(
(r_{\nu-1}-r_{\nu-2})\gamma_{\nu-1}\big) e(r_{\nu-1}\gamma_\nu).
\end{split}
\end{equation}
Equalities \eqref{2.6} and \eqref{2.3} further yield
\begin{equation*}
\begin{split}
\SSS_\nu & =\frac{1}{N} \sum\limits_{r=(r_1,\dots,r_{\nu-1})\in
\Z^{\nu -1}} \hspace{-10pt} c_r \hspace{-8pt}
\sum\limits_{\substack{d_1 \vert r_1 \\ d_2 \vert r_2-r_1 \\
\dots \\ d_{\nu-1}\vert r_{\nu-1}-r_{\nu-2} \\ d_\nu \vert
r_{\nu-1}}} \hspace{-10pt} d_1 \cdots d_{\nu} M\bigg(
\frac{Q}{d_1}\bigg)
\cdots M \bigg( \frac{Q}{d_{\nu}}\bigg) \\
& =\frac{1}{N} \sum\limits_{1\leq d_1,\dots,d_\nu\leq Q}
\hspace{-10pt} d_1\cdots d_\nu M\bigg( \frac{Q}{d_1}\bigg) \cdots
M\bigg( \frac{Q}{d_\nu}\bigg) \hspace{-10pt}
\sum\limits_{\substack{r\in
\Z^{\nu-1} \\ d_1 \vert r_1 \\ d_2 \vert r_2-r_1 \\
\dots \\ d_{\nu-1}\vert r_{\nu-1}-r_{\nu-2} \\ d_\nu \vert
r_{\nu-1}}} \hspace{-11pt} c_r .
\end{split}
\end{equation*}
The divisibility conditions $d_1 \vert r_1$, $d_2 \vert r_2-r_1$,
$d_3\vert r_3-r_2$,\dots,$d_{\nu-1}\vert r_{\nu-1}-r_{\nu-2}$,
$d_\nu\vert r_{\nu-1}$, read as
\begin{equation*}
\begin{split}
& r_1 =\ell_1 d_1 ,\\
& r_2=\ell_1 d_1+\ell_2 d_2 ,\\
& \dots \dots \dots \dots \dots \dots \dots \dots  \\
& r_{\nu-1} =\ell_1 d_1+\cdots +\ell_{\nu-1} d_{\nu-1} = \ell_\nu
d_\nu,
\end{split}
\end{equation*}
for some integers $\ell_1,\dots,\ell_\nu$. Thus, putting
$d=(d_1,\dots,d_{\nu-1})\in \Box_Q^{\nu-1} :=[1,Q]^{\nu-1}\cap
\Z^{\nu-1}$, $\ell =(\ell_1,\dots,\ell_{\nu-1})$, we can further
write
\begin{equation}\label{2.7}
\begin{split}
& \SSS_\nu =\frac{1}{N} \sum\limits_{d \in \Box_Q^{\nu-1}} d_1
\cdots d_{\nu-1} M\bigg( \frac{Q}{d_1}\bigg)
\cdots M\bigg( \frac{Q}{d_{\nu-1}} \bigg) \\
& \qquad \cdot \sum\limits_{\ell \in \Z^{\nu-1}} c_{(d_1\ell_1,d_1
\ell_1+d_2 \ell_2 , \dots,d_1\ell_1+\cdots
+d_{\nu-1}\ell_{\nu-1})} \hspace{-10pt} \sum\limits_{d_\nu \vert
d_1\ell_1+\cdots +d_{\nu-1}\ell_{\nu-1}} \hspace{-20pt} d_\nu
M\bigg( \frac{Q}{d_\nu}\bigg) .
\end{split}
\end{equation}

When $\nu=2$ we simply get
\begin{equation}\label{2.8}
\begin{split}
\SSS_2 & =\frac{1}{N} \sum\limits_{r\in \Z} c_r
\sum\limits_{\gamma_1,\gamma_2 \in \FF_Q} e\big(
r(\gamma_2-\gamma_1)\big) =\frac{1}{N} \sum\limits_{r\in \Z} c_r
\Bigg| \sum\limits_{\gamma \in \FF_Q} e(r\gamma)\bigg|^2 \\
& =\frac{1}{N} \sum\limits_{r\in \Z} c_r \Bigg|
\sum\limits_{\substack{1\leq d\leq Q \\ d\vert r}} dM\bigg(
\frac{Q}{d}\bigg) \Bigg|^2 \\ & =\frac{1}{N} \sum\limits_{r\in \Z}
c_r \sum\limits_{\substack{1\leq d_1,d_2 \leq Q \\ d_1 \vert r,\
d_2 \vert r}} d_1d_2 M\bigg( \frac{Q}{d_1}\bigg) M\bigg(
\frac{Q}{d_2}\bigg)\\
& =\frac{1}{N} \sum\limits_{1\leq d_1,d_2\leq Q} d_1 d_2 M\bigg(
\frac{Q}{d_1}\bigg) M\bigg( \frac{Q}{d_2}\bigg) \sum\limits_{\ell
\in \Z} c_{\ell [d_1,d_2]} .
\end{split}
\end{equation}

\section{Proof of Theorem \ref{T1}.}

In this section we will elaborate on formula \eqref{2.7} of the
smooth $\nu$-level correlation sum and prove that
\begin{equation}\label{3.1}
\lim\limits_{Q\rightarrow \infty}
R^{(\nu)}(Q,H)=2\sum\limits_{\substack{A_1,\dots,A_{\nu-1}\in \N \\
B_1,\dots,B_{\nu-1}\in \N  \\
(A_j,B_j)=1}} \ \ \ \iint\limits_{\Omega_{A,B,\Lambda}} H\big(
\Phi_{A,B}(x,y)\big) dx dy.
\end{equation}

Making use of \eqref{2.3} and of the fact that $\gamma \mapsto
1-\gamma$ is a bijection on $\FF_Q$, we can write the inner sum in
\eqref{2.7} as
\begin{equation*}
\sum\limits_{\gamma \in \FF_Q} e(\gamma d\cdot
\ell)=\sum\limits_{\gamma \in \FF_Q} e(-\gamma d\cdot \ell).
\end{equation*}
Taking also into account \eqref{2.5} we see that the contribution
of the two inner sums in \eqref{2.7} equals
\begin{equation*}
\begin{split}
& \sum\limits_{\substack{\ell \in \Z^{\nu-1} \\ \gamma \in \FF_Q}}
c_{(d_1\ell_1,d_1\ell_1+d_2\ell_2,\dots,d_1\ell_1+\cdots+d_{\nu-1}\ell_{\nu-1})}
 e(-\gamma d\cdot \ell) \\
& =\sum\limits_{\substack{\ell \in \Z^{\nu-1} \\ \gamma \in
\FF_Q}} \int\limits_{\R^{\nu-1}} e\Bigg(
-\gamma\sum\limits_{i=1}^{\nu -1} d_i \ell_i
-\sum\limits_{j=1}^{\nu-1} x_j (d_1\ell_1+\cdots+d_j
\ell_j)\Bigg)H(Nx) dx \\
& =\sum\limits_{\substack{\ell \in \Z^{\nu-1} \\ \gamma \in
\FF_Q}} \int\limits_{\R^{\nu-1}} e\Bigg(
-\sum\limits_{i=1}^{\nu-1} d_i\ell_i
(x_i+\cdots+x_{\nu-1}+\gamma))\Bigg)H(Nx) dx \\
& =\sum\limits_{\substack{\ell \in \Z^{\nu-1} \\ \gamma \in
\FF_Q}} \int\limits_{\R^{\nu-1}} e\Bigg(
-\sum\limits_{i=1}^{\nu-1} d_i\ell_i (x_i+\cdots+x_{\nu-1})\Bigg)
H\big( N(x_1,\dots,x_{\nu-2},x_{\nu-1}-\gamma) \big) dx.
\end{split}
\end{equation*}
Taking $y_i=d_i(x_i+\cdots+x_{\nu-1})$, $i=1,\dots ,\nu-1$, that
is
\begin{equation*}
\left\{ \begin{array}{l}
x_1 =\frac{y_1}{d_1}-\frac{y_2}{d_2} \\
x_2=\frac{y_2}{d_2}-\frac{y_3}{d_3} \\
\cdots \cdots \cdots \cdots \cdots \cdots  \\
x_{\nu-2}=\frac{y_{\nu-2}}{d_{\nu-2}} -\frac{y_{\nu-1}}{d_{\nu-1}} \\
x_{\nu-1} =\frac{y_{\nu-1}}{d_{\nu-1}} ,
\end{array}
\right.
\end{equation*}
and putting
\begin{equation*}
H_{N;d,\gamma} (y)= H\bigg( N\Big(
\frac{y_1}{d_1}-\frac{y_2}{d_2}\Big) , \dots ,N\Big(
\frac{y_{\nu-2}}{d_{\nu-2}}-\frac{y_{\nu-1}}{d_{\nu-1}}
\Big),N\Big( \frac{y_{\nu-1}}{d_{\nu-1}}-\gamma \Big) \bigg),
\end{equation*}
with $d=(d_1,\dots,d_{\nu-1})\in \Box_Q^{\nu-1}$,
$y=(y_1,\dots,y_{\nu-1})\in \R^{\nu-1}$, $\gamma \in \FF_Q$, the
contribution of the two inner sums in \eqref{2.7} becomes
\begin{equation*}
\begin{split}
\frac{1}{d_1\dots d_{\nu-1}} \sum\limits_{\gamma \in \FF_Q}
\sum\limits_{\ell \in \Z^{\nu-1}} \int\limits_{\R^{\nu-1}} &
e(-\ell \cdot y )\  H_{N;d,\gamma} (y) dy \\
& =\frac{1}{d_1 \dots d_{\nu-1}} \sum\limits_{\gamma \in \FF_Q}
\sum\limits_{\ell \in \Z^{\nu-1}} \widehat{H_{N;d,\gamma}} (\ell
).
\end{split}
\end{equation*}
Applying Poisson summation to the inner sum, this further equals
\begin{equation*}
\frac{1}{d_1\cdots d_{\nu-1}} \sum\limits_{\gamma \in \FF_Q}
\sum\limits_{\ell \in \Z^{\nu-1}} H_{N;d,\gamma} (\ell ),
\end{equation*}
which we insert back into \eqref{2.7} to get
\begin{equation*}
\SSS_\nu=\frac{1}{N} \sum\limits_{d \in \Box_Q^{\nu-1}} M\bigg(
\frac{Q}{d_1}\bigg) \cdots M\bigg( \frac{Q}{d_{\nu-1}} \bigg)
\sum\limits_{\gamma \in \FF_Q} \sum\limits_{\ell \in \Z^{\nu-1}}
H_{N;d ,\gamma} (\ell).
\end{equation*}
The support of $H$ is included in $(0,\Lambda^\prime)^{\nu-1}$,
thus we necessarily have
\begin{equation*}
0<N\bigg( \frac{\ell_j}{d_j}-\frac{\ell_{j+1}}{d_{j+1}} \bigg)
<\Lambda^\prime ,\qquad j=1,\dots, \nu-2.
\end{equation*}
These inequalities firstly give that $\ell_j d_{j+1}-\ell_{j+1}d_j
\geq 1$, and secondly, that
\begin{equation*}
\Lambda^\prime >\frac{N(\ell_j d_{j+1}-\ell_{j+1}d_j)}{d_jd_{j+1}}
\geq \frac{N}{d_jd_{j+1}} .
\end{equation*}
Therefore for all $Q\geq Q_0(\Lambda^\prime)$ we find that
\begin{equation*}
\frac{Q^2}{d_jd_{j+1}} =\frac{Q^2}{N}\cdot \frac{N}{d_jd_{j+1}}
<\frac{Q^2 \Lambda^\prime}{N} <c_\Lambda=\frac{\pi^2 \Lambda}{3} .
\end{equation*}
Here both $\frac{Q}{d_j}$ and $\frac{Q}{d_{j+1}}$ are $\geq 1$.
Hence each of them is $\leq c_\Lambda$. It follows that for $Q\geq
Q_0(\Lambda^\prime )$ we have
\begin{equation}\label{3.2}
1\leq \frac{Q}{d_j} \leq c_\Lambda , \qquad j=1,\dots,\nu-1.
\end{equation}
Same reasoning for $j=\nu-1$ gives, for $Q\geq
Q_0(\Lambda^\prime)$,
\begin{equation}\label{3.3}
\frac{Q}{q} \leq c_\Lambda .
\end{equation}
Therefore
\begin{equation*}
\begin{split}
\SSS_\nu & =\frac{1}{N} \sum\limits_{d \in \Box_Q^{\nu-1}}
\sum\limits_{\ell \in \Z^{\nu-1}} \sum\limits_{\substack{1\leq
r_j\leq \frac{Q}{d_j}\\ j=1,\dots,\nu-1}} \mu (r_1)\cdots
\mu(r_{\nu-1}) \sum\limits_{\substack{\frac{a}{q}\in \FF_Q \\
q\geq \frac{Q}{c_\Lambda}}}
H_{N;d ,\frac{a}{q}} (\ell)   \\
& =\frac{1}{N}\sum\limits_{1\leq r_1,\dots,r_{\nu-1}\leq
c_\Lambda} \mu(r_1)\cdots \mu(r_{\nu-1})
\sum\limits_{\substack{1\leq d_j \leq \frac{Q}{r_j}\\
j=1,\dots,\nu-1}} \ \sum\limits_{\ell \in \Z^{\nu-1}} \
\sum\limits_{\substack{\frac{a}{q}\in \FF_Q \\
q\geq \frac{Q}{c_\Lambda}}} H_{N;d ,\frac{a}{q}} (\ell).
\end{split}
\end{equation*}
The inner sum above is given by
\begin{equation*}
\sum\limits_{\substack{\frac{a}{q}\in \FF_Q \\ q\geq
\frac{Q}{c_\Lambda} }} H\Bigg( N\bigg(
\frac{\ell_1}{d_1}-\frac{\ell_2}{d_2}\bigg) ,\dots, N\bigg(
\frac{\ell_{\nu-2}}{d_{\nu-2}}-\frac{\ell_{\nu-1}}{ d_{\nu-1}}
\bigg),N\bigg( \frac{\ell_{\nu-1}}{d_{\nu-1}}- \frac{a}{q}\bigg)
\Bigg) .
\end{equation*}

For $j=1,\dots,\nu-1$, we set
\begin{equation}\label{3.4}
\Delta_j=q\ell_j-ad_j .
\end{equation}
Since $\supp (H) \subset (0,\Lambda^\prime]^{\nu-1}$, we have
\begin{equation*}
0<\frac{N\Delta_j}{qd_j}=N\bigg( \frac{\ell_j}{d_j}-\frac{a}{q}
\bigg) =N\bigg( \frac{\ell_j}{d_j}-\frac{\ell_{j+1}}{d_{j+1}}
\bigg)+\cdots +N\bigg( \frac{\ell_{\nu-1}}{d_{\nu-1}}-\frac{a}{q}
\bigg) <(\nu-j)\Lambda^\prime .
\end{equation*}
Hence $\Delta_j \geq 1$ and, for $Q\geq Q_0(\Lambda^\prime)$,
\begin{equation*}
\Delta_j \leq \frac{qd_j (\nu-j)\Lambda^\prime}{N} \leq
\frac{Q^2(\nu-j)\Lambda^\prime}{N} \leq (\nu-j)c_\Lambda ,
\end{equation*}
therefore
\begin{equation*}
1\leq \Delta_1,\dots,\Delta_{\nu-1} \leq (\nu-1)c_\Lambda .
\end{equation*}
Note also from \eqref{3.4} that $\ell_j$ is uniquely determined as
\begin{equation*}
\ell_j=\frac{\Delta_j+ad_j}{q} .
\end{equation*}
This gives in turn that
\begin{equation*}
\frac{\ell_j}{d_j}-\frac{\ell_{j+1}}{d_{j+1}}
=\frac{\Delta_j+ad_j}{qd_j}-
\frac{\Delta_{j+1}+ad_{j+1}}{qd_{j+1}} = \frac{1}{q}\bigg(
\frac{\Delta_j}{d_j}- \frac{\Delta_{j+1}}{d_{j+1}} \bigg), \quad
j=1,\dots,\nu-2.
\end{equation*}
We also have
\begin{equation*}
\frac{\ell_{\nu-1}}{d_{\nu-1}}-\frac{a}{q} =
\frac{\Delta_{\nu-1}}{qd_{\nu-1}} .
\end{equation*}
Here $d_j$ needs to satisfy the congruence
\begin{equation*}
d_j =-\bar{a} \Delta_j \pmod{q} ,\qquad j=1,\dots,\nu-1,
\end{equation*}
where $\bar{a}$ denotes the integer between $1$ and $q$ which
satisfies $a\bar{a}=1\pmod{q}$.

In summary, we infer that
\begin{equation}\label{3.5}
\begin{split}
\SSS_\nu =\frac{1}{N} & \sum\limits_{1\leq r_1,\dots,r_{\nu-1}
\leq c_\Lambda} \mu(r_1)\cdots \mu(r_{\nu-1}) \sum\limits_{1\leq
\Delta_1,\dots,\Delta_{\nu-1} \leq (\nu -1)c_\Lambda}
\sum\limits_{\substack{\frac{a}{q}\in \FF_Q \\
q\geq \frac{Q}{c_\Lambda} }} \\ & \sum\limits_{\substack{1\leq
d_j\leq \frac{Q}{r_j} \\ d_j=-\bar{a}\Delta_j
\hspace{-8pt}\pmod{q} \\j=1,\dots,\nu-1}} \hspace{-20pt} H\Bigg(
\frac{N}{q} \bigg(
\frac{\Delta_1}{d_1}-\frac{\Delta_2}{d_2},\dots,
\frac{\Delta_{\nu-2}}{d_{\nu-2}}-
\frac{\Delta_{\nu-1}}{d_{\nu-1}},
\frac{\Delta_{\nu-1}}{d_{\nu-1}}\bigg) \Bigg) .
\end{split}
\end{equation}
To simplify this expression, we consider the linear transformation
$T$ defined by \eqref{1.3}. The function $\tilde{H}=H\circ T$ is
smooth, $\supp (\tilde{H}) \subset \big(0,(\nu-1)\Lambda^\prime ]
\times \cdots \times (0,\Lambda^\prime ]$, and the inner sum in
\eqref{3.5} becomes
\begin{equation*}
\sum\limits_{\substack{1\leq d_j\leq \frac{Q}{r_j} \\
d_j=-\bar{a}\Delta_j \hspace{-8pt}\pmod{q} \\ j=1,\dots,\nu-1}}
\hspace{-5pt} \tilde{H}\Bigg( \frac{N}{q} \bigg(
\frac{\Delta_1}{d_1}, \frac{\Delta_{2}}{d_{2}},
\dots,\frac{\Delta_{\nu-1}}{d_{\nu-1}}\bigg) \Bigg) .
\end{equation*}

When $c_\Lambda <1$, ${\mathcal S}_\nu=0$. For $j=1,\dots,\nu-1$,
we define
\begin{equation*}
e_j=\frac{d_j+\bar{a}\Delta_j}{q} .
\end{equation*}
The congruence $d_j=-\bar{a}\Delta_j \hspace{-3pt} \pmod{q}$ shows
that each $e_j$ is an integer. Moreover, $d_j,\bar{a},\Delta_j$
are all greater or equal than $1$, so $e_j \geq 1$. On the other
hand, using \eqref{3.3}, we see that $\frac{d_j}{q}\leq
\frac{Q}{qr_j}\leq \frac{Q}{q}\leq c_\Lambda$ and
$\frac{\bar{a}\Delta_j}{q}<\Delta_j \leq (\nu-1)c_\Lambda$,
leading to
\begin{equation*}
1\leq e_j \leq \nu c_\Lambda ,\qquad j=1,\dots,\nu -1.
\end{equation*}
With $q$, $a$, $\Delta_j$ fixed, each value of $e_j$ uniquely
determines a value of $d_j$, precisely
\begin{equation*}
d_j=qe_j -\bar{a} \Delta_j .
\end{equation*}
Moreover, with $e_j$ and $\Delta_j$ fixed and $\frac{a}{q}$
variable in $\FF_Q$, $a$ and $q$ need to satisfy some extra
conditions in order for $d_j$ to belong to the set $\big\{
1,\dots,\big[ \frac{Q}{r_j}\big]\big\}$. Using \eqref{3.2}, we
infer that $a$ and $q$ necessarily fulfil
\begin{equation*}
\frac{Q}{c_\Lambda r_j}\leq \frac{Q}{c_\Lambda} \leq qe_j-\bar{a}
\Delta_j \leq \frac{Q}{r_j},\qquad j=1,\dots,\nu -1.
\end{equation*}

Consider now the convex region $\Omega_{r,e,\Delta}$ in $\R^2$
defined by the inequalities
\begin{equation*}
\begin{split}
& 0\leq x\leq y\leq 1, \qquad y\geq \frac{1}{c_\Lambda},\\ &
\frac{1}{c_\Lambda r_j} \leq e_j y-\Delta_j x \leq
\frac{1}{r_j},\quad j=1,\dots,\nu-1 ,
\end{split}
\end{equation*}
and the functions $f_{e,\Delta},f_{e,\Delta}^{(j)}$ defined on
$Q\Omega_{r,e,\Delta}$ by
\begin{equation*}
\begin{split}
& f_{e,\Delta}(b,q) =\tilde{H}\Big( f^{(1)}_{e,\Delta}
(b,q),\dots,f^{(\nu -1)}_{e,\Delta}(b,q)\Big),\\
& f^{(j)}_{e,\Delta}(b,q) =\frac{N\Delta_j}{q(qe_j-b\Delta_j)} ,
\qquad j=1,\dots, \nu -1.
\end{split}
\end{equation*}
We write $b=\bar{a}$, and remark that as $\frac{a}{q}$ runs over
$\FF_Q$ with $q\geq \frac{Q}{c_\Lambda}$, so does $\frac{b}{q}$.
Thus
\begin{equation*}
\SSS_\nu =\frac{1}{N} \sum\limits_{1\leq r_1,\dots,r_{\nu-1} \leq
c_\Lambda} \mu(r_1) \cdots \mu(r_{\nu-1})
\hspace{-6pt}\sum\limits_{\substack{1\leq
\Delta_1,\dots,\Delta_{\nu-1}\leq (\nu -1)c_\Lambda \\ 1\leq
e_1,\dots,e_{\nu-1} \leq \nu c_\Lambda}}
\sum\limits_{\substack{(b,q)\in Q\Omega_{r ,e,\Delta} \\
(b,q)=1}} \hspace{-6pt} f_{e,\Delta}(b,q).
\end{equation*}
By Corollary 1 in \cite{BCZ} the inner sum above can be written as
\begin{equation}\label{3.6}
\frac{6}{\pi^2} \iint\limits_{Q\Omega_{r,e,\Delta}}
f_{e,\Delta}(u,v) du dv\ +O\bigg( \frac{\| Df_{e,\Delta}\|_\infty
Q^2 \log Q+\| f_{e,\Delta}\|_\infty Q\log Q}{N}\bigg).
\end{equation}
Since $H=\tilde{H}\circ T^{-1}$, it is clear that $\|
f_{e,\Delta}\|_\infty \leq \| H\|_\infty$. Using the definition of
$\Omega_{r,e,\Delta}$, we also find, for every $j=1,\dots, \nu
-1$, that
\begin{equation*}
\begin{split}
\| Df_{e,\Delta} \|_\infty & =\sup\limits_{(b,q)\in
Q\Omega_{r,e,\Delta}} \bigg( \ \bigg| \frac{\partial
f^{(j)}_{e,\Delta}}{\partial b}\bigg| +\bigg| \frac{\partial
f^{(j)}_{e,\Delta}}{\partial q}\bigg|\ \bigg) \ll_\Lambda
N\sup\limits_{(b,q)\in Q\Omega_{r,e,\Delta}} \frac{1}{q(qe_j-b\Delta_j)^2}\\
& \ll \frac{1}{Q} \sup\limits_{(x,y)\in \Omega_{r,e,\Delta}}
\frac{1}{y(e_jy-\Delta_j x)^2} \ll_\Lambda \frac{1}{Q} ,
\end{split}
\end{equation*}
showing that the error term in \eqref{3.6} is $\ll_H \frac{\log
Q}{Q}$. Rescaling to $(u,v)=(Qx,Qy)$ we find that
\begin{equation}\label{3.7}
\SSS_\nu =\frac{6Q^2}{\pi^2 N} \sum\limits_{1\leq
r_1,\dots,r_{\nu-1} \leq c_\Lambda} \hspace{-6pt} \mu(r_1) \cdots
\mu(r_{\nu -1}) \hspace{-6pt}\sum\limits_{\substack{1\leq
\Delta_j\leq (\nu -1)c_\Lambda \\ 1\leq e_j \leq \nu c_\Lambda
\\ j=1,\dots,\nu-1}} \hspace{-6pt} I_{r,e,\Delta}+O_H \bigg(
\frac{\log Q}{Q}\bigg) ,
\end{equation}
where this time we put
\begin{equation}\label{3.8}
\begin{split}
I_{r,e,\Delta} & =\iint\limits_{\Omega_{r,e,\Delta}}
g_{e,\Delta}(x,y)dx dy,
\\ g_{e,\Delta}(x,y) & =\tilde{H}\Big(
g^{(1)}_{e,\Delta}(x,y),\dots,g^{(\nu-1)}_{e,\Delta}(x,y)\Big) ,\\
g^{(j)}_{e,\Delta}(x,y) & =\frac{N\Delta_j}{Q^2 y(e_j y-\Delta_j
x)} ,\qquad j=1,\dots,\nu-1.
\end{split}
\end{equation}
Using \eqref{2.1} and the inequality
\begin{equation*}
\big| \tilde{H} (v )-\tilde{H}(w) \big| \leq \| \tilde{H}^\prime
\| \, \vert v-w\vert \leq 2\| H^\prime \|  \,\vert v-w \vert,
\end{equation*}
we see that formula \eqref{3.7} holds true after replacing
$g^{(j)}_{e,\Delta}$ by
\begin{equation}\label{3.9}
g^{(j)}_{e,\Delta}(x,y) =\frac{3\Delta_j}{\pi^2 y(e_jy-\Delta_j
x)} ,\qquad j=1,\dots, \nu-1,
\end{equation}
in the formula for $g_{e,\Delta}$ from \eqref{3.8}. Therefore we
infer that
\begin{equation}\label{3.10}
\SSS_\nu=2\sum\limits_{1\leq r_1,\dots,r_{\nu-1} \leq c_\Lambda}
\hspace{-6pt} \mu(r_1) \cdots \mu(r_{\nu-1})
\hspace{-6pt}\sum\limits_{\substack{1\leq
\Delta_1,\dots,\Delta_{\nu-1}\leq (\nu-1)c_\Lambda \\ 1\leq
e_1,\dots,e_{\nu-1} \leq \nu c_\Lambda}} \hspace{-6pt}
I_{r,e,\Delta}+O_H \bigg( \frac{\log Q}{Q}\bigg) ,
\end{equation}
where $I_{r,e,\Delta}$ and $g_{e,\Delta}$ are as in \eqref{3.8}
and $g^{(j)}_{e,\Delta}$ as in \eqref{3.9}.

Next, we notice that the region $\Omega_{r,e,\Delta}$ can be
extended to
\begin{equation*}
\tilde{\Omega}_{r,e,\Delta} =\bigg\{ (x,y) : 0\leq x\leq y\leq 1,
y\geq \frac{1}{c_\Lambda},  0<e_j y-\Delta_j x\leq
\frac{1}{r_j},j=1,\dots,\nu-1 \bigg\} ,
\end{equation*}
without changing the terms $I_{r,e,\Delta}$ in \eqref{3.10}.
Indeed, if $(x,y)\in \tilde{\Omega}_{r,e,\Delta} \setminus
\Omega_{r,e,\Delta}$, there is $j$ for which $\vert e_j y-\Delta_j
x\vert <\frac{1}{c_\Lambda r_j}$, and thus
\begin{equation*}
\vert g_{e,\Delta}^{(j)} (x,y)\vert \geq \frac{3\Delta_j}{
\frac{\pi^2}{c_\Lambda r_j}} =\frac{3c_\Lambda \Delta_j
r_j}{\pi^2} \geq \frac{3c_\Lambda}{\pi^2}=\Lambda .
\end{equation*}
This yields that $g_{e,\Delta}=0$ on $\tilde{\Omega}_{r,e,\Delta}
\setminus \Omega_{r,e,\Delta}$. Hence
\begin{equation}\label{3.11}
\begin{split}
\SSS_\nu & =2 \sum\limits_{1\leq r_1,\dots,r_{\nu-1} \leq
c_\Lambda} \hspace{-6pt} \mu(r_1) \cdots \mu(r_{\nu-1}) \\ &
\qquad \sum\limits_{\substack{1\leq
\Delta_1,\dots,\Delta_{\nu-1}\leq (\nu-1)c_\Lambda \\ 1\leq
e_1,\dots,e_{\nu-1} \leq \nu c_\Lambda}}\ \
\iint\limits_{\tilde{\Omega}_{r,e,\Delta}} g_{e,\Delta}(x,y) dx dy
+O_H \bigg( \frac{\log Q}{Q}\bigg) .
\end{split}
\end{equation}

We take $A_j=e_j r_j$, $B_j=\Delta_j r_j$,
$A=(A_1,\dots,A_{\nu-1})$, $B=(B_1,\dots,B_{\nu-1})$, and consider
regions $\Omega_{A,B,\Lambda}$ and maps $T_{A,B}$ as defined in
\eqref{1.2} and \eqref{1.1}. We have
$\Omega_{A,B,\Lambda}=\tilde{\Omega}_{r,e,\Delta}$. If we denote
\begin{equation*}
I_{A,B,\Lambda} =\iint\limits_{\Omega_{A,B,\Lambda}}
\tilde{H}\big( T_{A,B}(x,y)\big) dx
dy=\iint\limits_{\Omega_{A,B,\Lambda}} H\big( \Phi_{A,B}(x,y)\big)
dx dy ,
\end{equation*}
then \eqref{3.11} yields that
\begin{equation*}
\begin{split}
\SSS_\nu & =2\sum\limits_{\substack{1\leq A_1,\dots,A_{\nu-1}\leq
(\nu-1)c_\Lambda^2 \\ 1\leq B_1,\dots,B_{\nu-1} \leq \nu
c_\Lambda^2}} \hspace{-6pt}
I_{A,B,\Lambda} \sum\limits_{\substack{r_j \vert (A_j,B_j) \\
j=1,\dots,\nu-1}} \mu(r_1)\dots \mu(r_{\nu-1})
 +O_H \bigg( \frac{\log Q}{Q}\bigg) \\ & =2
 \sum\limits_{\substack{1\leq A_1,\dots,A_{\nu-1}\leq
 (\nu-1)c_\Lambda^2 \\ 1\leq B_1,\dots,B_{\nu-1}\leq \nu
c_\Lambda^2 \\ (A_j,B_j)=1,j=1,\dots,\nu-1}} \hspace{-6pt}
I_{A,B,\Lambda}+O_H \bigg( \frac{\log Q}{Q}\bigg) \\
& =2\sum\limits_{\substack{A,B\in \N^{\nu-1} \\
(A_j,B_j)=1}} I_{A,B,\Lambda}+O_H \bigg( \frac{\log Q}{Q}\bigg),
\end{split}
\end{equation*}
which concludes the proof of \eqref{3.1}.

Theorem \ref{T1} now follows from \eqref{3.1} by approximating
pointwise the characteristic function $\chi_\BB$ of a bounded box
$\BB \subset (0,\infty)^{\nu-1}$ by smooth functions with compact
support $H_\pm$ such that $\supp (H_\pm) \subset
(0,\Lambda)^{\nu-1}$ and $0\leq H_- \leq \chi_\BB \leq H_+ \leq
1$.

\section{Proof of Theorem \ref{T2}.}

To establish the formula for $g_2$ given in \eqref{1.5}, we return
to the formula for $\SSS_2$ from \eqref{2.8} and consider for each
$y>0$ the function
\begin{equation*}
H_y(x)=\frac{1}{y} H\bigg( \frac{x}{y}\bigg),\qquad x\in \R .
\end{equation*}
Then
\begin{equation*}
\widehat{H_y} (z)=\widehat{H}(yz),
\end{equation*}
and from \eqref{2.5} we find that the inner sum in \eqref{2.8} can
be written as
\begin{equation}\label{4.1}
\frac{1}{N} \sum\limits_{\ell \in \Z} \widehat{H} \bigg(
\frac{\ell [d_1, d_2]}{N} \bigg)=\frac{1}{N} \sum\limits_{\ell \in
\Z} \widehat{H_{\frac{[d_1,d_2]}{N}}} (\ell).
\end{equation}
By Poisson's summation formula we have
\begin{equation}\label{4.2}
\sum\limits_{\ell \in \Z} \widehat{H_{\frac{[d_1,d_2]}{N}}} (\ell)
=\sum\limits_{\ell \in \Z} H_{\frac{[d_1,d_2]}{N}} (\ell).
\end{equation}
Combining \eqref{2.8}, \eqref{4.1} and \eqref{4.2} we find that
\begin{equation}\label{4.3}
\SSS_2=\frac{1}{N^2} \sum\limits_{1\leq d_1,d_2\leq Q} d_1d_2
M\bigg( \frac{Q}{d_1}\bigg)\, M\bigg( \frac{Q}{d_2}\bigg)
\sum\limits_{\ell \in \Z} H_{\frac{[d_1,d_2]}{N}} (\ell).
\end{equation}

Using the definition of $M$ and $H_y$ we can rewrite \eqref{4.3}
as
\begin{equation}\label{4.4}
\begin{split}
\SSS_2 & =\frac{1}{N^2} \sum\limits_{1\leq d_1,d_2 \leq Q}
\hspace{-6pt} d_1 d_2
\sum\limits_{\substack{1\leq r_1 \leq \frac{Q}{d_1} \\
1\leq r_2 \leq \frac{Q}{d_2}}} \mu(r_1)\mu (r_2) \sum\limits_{\ell
\in \Z} \frac{N}{[d_1,d_2]} H\bigg(
\frac{\ell N}{[d_1,d_2]}\bigg) \\
& =\frac{1}{N} \sum\limits_{\substack{1\leq r_1d_1 \leq Q \\
1\leq r_2d_2 \leq Q}} \mu(r_1)\mu(r_2) (d_1,d_2) \sum\limits_{\ell
\in \Z} H\bigg( \frac{\ell N}{[d_1,d_2]}\bigg)\\
& =\frac{1}{N} \sum\limits_{\substack{1\leq r_1\leq Q \\ 1\leq r_2
\leq Q}} \mu(r_1)\mu(r_2) \sum\limits_{\substack{1\leq d_1 \leq
\frac{Q}{r_1} \\ 1\leq d_2 \leq \frac{Q}{r_2}}} (d_1,d_2)
\sum\limits_{\ell \in \Z} H\bigg( \frac{\ell N}{[d_1,d_2]}\bigg) .
\end{split}
\end{equation}
Denote $\delta=(d_1,d_2)$, so that $d_1=\delta q_1$, $d_2=\delta
q_2$ with $(q_1,q_2)=1$. Then \eqref{4.4} becomes
\begin{equation}\label{4.5}
\SSS_2 =\frac{1}{N} \sum\limits_{1\leq r_1,r_2\leq Q}\
\sum\limits_{1\leq \delta \leq \frac{Q}{\max \{ r_1,r_2\}}}
\hspace{-6pt} \mu(r_1)\mu(r_2)\delta \hspace{-6pt}
\sum\limits_{\substack{1\leq q_1\leq \frac{Q}{\delta r_1}
\\ 1\leq q_2 \leq \frac{Q}{\delta r_2} \\ (q_1,q_2)=1}}
\sum\limits_{\ell \in \Z} H\bigg( \frac{\ell N}{\delta q_1
q_2}\bigg) .
\end{equation}

Only values of $\ell$ with $\ell <\frac{\delta q_1q_2
\Lambda^\prime}{N}$ may produce a non-zero contribution in the
inner sum in \eqref{4.5}. For $Q$ larger than some
$Q_0(\Lambda^\prime)$ this yields in conjunction with \eqref{2.1}
that
\begin{equation*}
\delta r_1r_2 \ell <\frac{\delta^2 r_1r_2q_1q_2\Lambda^\prime}{N}
\leq \frac{Q^2\Lambda^\prime}{N}<c_\Lambda =\frac{\pi^2
\Lambda}{3} .
\end{equation*}
So each of $\delta,r_1,r_2,\ell$ should be at most $c_\Lambda$ and
in \eqref{4.5} we are left with
\begin{equation*}
\SSS_2=\frac{1}{N} \sum\limits_{\substack{\ell,\delta,r_1,r_2 \geq 1\\
\ell \delta r_1r_2 <c_\Lambda}} \hspace{-6pt} \mu(r_1)\mu(r_2)
\delta \sum\limits_{\substack{1\leq q_1\leq \frac{Q}{\delta r_1} \\
1\leq q_2 \leq \frac{Q}{\delta r_2} \\ (q_1,q_2)=1}} H\bigg(
\frac{\ell N}{\delta q_1q_2}\bigg) .
\end{equation*}
Here
\begin{equation*}
H\bigg( \frac{\ell N}{\delta q_1q_2}\bigg)=H\bigg( \frac{3\ell
Q^2}{\pi^2 \delta q_1q_2}\bigg)+O_\Lambda (\| H^\prime \|_\infty
Q\log Q),
\end{equation*}
thus
\begin{equation*}
\SSS_2=\frac{1}{N} \sum\limits_{\substack{\ell,\delta,r_1,r_2 \geq 1\\
\ell \delta r_1r_2 <c_\Lambda}} \hspace{-6pt} \mu(r_1)\mu(r_2)
\delta \hspace{-6pt}
\sum\limits_{\substack{1\leq q_1\leq \frac{Q}{\delta r_1} \\
1\leq q_2 \leq \frac{Q}{\delta r_2} \\ (q_1,q_2)=1}}
\hspace{-10pt} H\bigg( \frac{3\ell Q^2}{\pi^2 \delta q_1q_2}\bigg)
+O_H \bigg( \frac{\log Q}{Q}\bigg) .
\end{equation*}
If $\min \{ q_1,q_2 \} <Q^{1-\eps}$ for some $\eps >0$, then for
sufficiently large $Q$ (with respect to $\eps$ and $\Lambda$) one
has $\frac{3\ell Q^2}{\pi^2 \delta q_1 q_2}>\Lambda$, and thus
$H\big( \frac{3\ell Q^2}{\pi^2 \delta q_1q_2}\big)=0$. On the
other hand we have (see Corollary 1 in \cite{BCZ})
\begin{equation}\label{4.6}
\sum\limits_{\substack{\min \{ q_1,q_2\} \geq Q^{1-\eps} \\
q_1 \leq \frac{Q}{\delta r_1} ,\, q_2 \leq \frac{Q}{\delta r_2}\\
(q_1,q_2)=1}} \hspace{-6pt} H\bigg( \frac{3\ell Q^2}{\pi^2 \delta
q_1q_2}\bigg) = \frac{6}{\pi^2} \hspace{-6pt}
\iint\limits_{\substack{x\leq \frac{Q}{\delta r_1} \\
y\leq \frac{Q}{\delta r_2} \\ \min \{ x,y\} \geq Q^{1-\eps}}}
\hspace{-10pt} H\bigg( \frac{3\ell Q^2}{\pi^2 \delta xy}\bigg) dx
dy +E_{H,\eps} (Q),
\end{equation}
where
\begin{equation*}
E_{H,\eps} (Q) \ll_{H,\eps} \frac{Q^2 \|
H^\prime\|_\infty}{Q^{3(1-\eps)}} Q^2 \log Q+\| H\|_\infty Q\log Q
\ll_H Q^{1+4\eps} \log Q.
\end{equation*}
The change of variables $(x,y)=(Qu,Qv)$ gives that the main term
in \eqref{4.6} can be expressed as
\begin{equation}\label{4.7}
\frac{6Q^2}{\pi^2} \hspace{-6pt}
\iint\limits_{\substack{u\leq \frac{1}{\delta r_1} \\
v\leq \frac{1}{\delta r_2} \\ \min \{ u,v\} \geq Q^{-\eps}}}
\hspace{-6pt} H\bigg( \frac{3\ell}{\pi^2 \delta uv}\bigg) du dv.
\end{equation}
For $u,v$ as in \eqref{4.7} and $Q\geq Q_0(\Lambda)$ we have
$\frac{\pi^2 \delta uv}{3\ell} \leq \frac{\pi^2}{3}Q^{-\eps}
<\frac{1}{\Lambda}$ (since $r_1,r_2,\ell \geq 1$); thus
\begin{equation}\label{4.8}
\iint\limits_{\substack{u\leq \frac{1}{\delta r_1} \\
v\leq \frac{1}{\delta r_2} \\ \min \{ u,v\} \leq Q^{-\eps}}}
\hspace{-12pt} H\bigg( \frac{3\ell}{\pi^2 \delta uv}\bigg) du
dv=0.
\end{equation}
As a result of \eqref{4.6}--\eqref{4.8} and of \eqref{2.1} we
infer that
\begin{equation}\label{4.9}
\SSS_2=2\hspace{-8pt}
\sum\limits_{\substack{\ell, \delta,r_1,r_2 \geq 1 \\
\ell \delta r_1r_2 <c_\Lambda}} \hspace{-6pt} \mu(r_1)\mu(r_2)
\delta \int\limits_0^{\frac{1}{\delta r_1}}
\int\limits_0^{\frac{1}{\delta r_2}} H\bigg( \frac{3\ell}{\pi^2
\delta uv}\bigg) dv du+O_{H,\eps} ( Q^{-1+\eps}\log Q).
\end{equation}
Next, we put $\lambda=\lambda_u(v)=\frac{3\ell}{\pi^2 \delta uv}$
and change the order of integration to express the double integral
in \eqref{4.9} as
\begin{equation*}
\begin{split}
\int\limits_0^{\frac{1}{\delta r_1}} \int\limits_{\frac{3\ell
r_2}{\pi^2 u}}^\Lambda H(\lambda)\frac{3\ell}{\pi^2 \delta
u\lambda^2} \ d\lambda du & =\frac{3\ell}{\pi^2 \delta}
\int\limits_\frac{3\ell r_2\delta r_1}{\pi^2}^\Lambda
\int\limits_\frac{3\ell r_2}{\pi^2 \lambda}^\frac{1}{\delta r_1}
\frac{H(\lambda)}{\lambda^2 u}\ du d\lambda \\
& =\frac{3\ell}{\pi^2 \delta} \int\limits_\frac{3\ell r_1r_2
\delta}{\pi^2}^\Lambda \frac{H(\lambda)}{\lambda^2} \log
\frac{\pi^2 \lambda}{3\ell r_1r_2\delta} d\lambda .
\end{split}
\end{equation*}
Inserting this back into \eqref{4.9} we get
\begin{equation}\label{4.10}
\SSS_2=\frac{6}{\pi^2}
\sum\limits_{\substack{\ell, \delta,r_1,r_2\geq 1 \\
\ell \delta r_1r_2 <c_\Lambda}} \hspace{-6pt} \mu(r_1)\mu(r_2)\ell
\hspace{-6pt} \int\limits_\frac{3\ell r_1r_2\delta}{\pi^2}^\Lambda
\hspace{-6pt} \frac{H(\lambda)}{\lambda^2} \log \frac{\pi^2
\lambda}{3\ell r_1r_2\delta} d\lambda +O_{H,\eps} ( Q^{-1+\eps}).
\end{equation}
At this point we take $K=\ell \delta r_1r_2$ and rewrite
\eqref{4.10} as
\begin{equation}\label{4.11}
\SSS_2=\frac{6}{\pi^2} \sum\limits_{1\leq K<c_\Lambda}\
\int\limits_{\frac{3K}{\pi^2}}^\Lambda
\frac{H(\lambda)}{\lambda^2}\log \frac{\pi^2 \lambda}{3K} d\lambda
\hspace{-6pt} \sum\limits_{\substack{\ell, \delta ,r_1,r_2 \geq 1 \\
\ell \delta r_1 r_2 =K}} \hspace{-16pt} \mu (r_1)\mu(r_2) \ell +
O_{H,\eps} ( Q^{-1+\eps} ) .
\end{equation}
Next, we put the inner sum on the right-hand side of \eqref{4.11}
in the form
\begin{equation}\label{4.12}
\begin{split}
\sum\limits_{\substack{\ell,r_1,r_2 \geq 1 \\
\ell r_1 r_2 \vert K}} \mu(r_1)\mu(r_2) \ell & =
\sum\limits_{\substack{\ell, r_1 \geq 1 \\
\ell r_1 \vert K}} \mu(r_1) \ell \sum\limits_{r_2 \vert
\frac{K}{\ell r_1}} \mu (r_2) =\sum\limits_{\substack{\ell ,r_1 \geq 1 \\
\ell r_1 = K}} \mu(r_1) \ell \\ & =K\sum\limits_{r_1 \vert K}
\frac{\mu(r_1)}{r_1} =K\frac{\varphi (K)}{K}=\varphi (K).
\end{split}
\end{equation}
Combining \eqref{4.11} and \eqref{4.12} we find that
\begin{equation*}
\begin{split}
\SSS_2 & =R^{(2)}(Q,H)=\frac{6}{\pi^2} \sum\limits_{1\leq
K<c_\Lambda} \varphi (K) \int\limits_{\frac{3K}{\pi^2}}^\Lambda
\frac{H(\lambda)}{\lambda^2} \log \frac{\pi^2 \lambda}{3K}
d\lambda +O_{H,\eps} ( Q^{-1+\eps}) \\ & =\frac{6}{\pi^2}
\int\limits_0^\Lambda \frac{H(\lambda)}{\lambda^2}
\sum\limits_{1\leq K<\frac{\pi^2 \Lambda}{3}} \varphi (K) \max
\bigg\{ 0,\log \frac{\pi^2 \lambda}{3K}\bigg\} d\lambda
+O_{H,\eps} ( Q^{-1+\eps})
\\ & =\int\limits_0^\Lambda H(\lambda)g_2 (\lambda)
d\lambda+O_{H,\eps}(Q^{-1+\eps}).
\end{split}
\end{equation*}
We now conclude the proof of \eqref{1.5} by approximating
pointwise the characteristic function of any closed interval
$I\subset (0,\Lambda)$ from below and from above by smooth
functions $H_\pm$ with support in $(0,\Lambda)$ and such that
$0\leq H_- \leq \chi_I \leq H_+ \leq 1$.

To prove \eqref{1.6} we start with the well-known equalities
\begin{equation}\label{4.13}
\sum\limits_{q=1}^\infty \frac{\varphi(q)}{q^s}
=\frac{\zeta(s-1)}{\zeta(s)}\qquad (\Real s>2),
\end{equation}
and
\begin{equation}\label{4.14}
\frac{1}{2\pi i} \int\limits_{\sigma_0-i\infty}^{\sigma_0+i\infty}
\frac{y^s}{s^2} ds =\begin{cases} 0, & 0\leq y\leq 1\\
\log y, & y>1\end{cases} \qquad (\sigma_0 >2).
\end{equation}
From \eqref{4.13} and \eqref{4.14} we infer for fixed $\sigma_0
>2$ that
\begin{equation}\label{4.15}
\sum\limits_{1\leq q<x} \varphi (q) \log \frac{x}{q}=\frac{1}{2\pi
i} \int\limits_{\sigma_0-i\infty}^{\sigma_0+i\infty}
\frac{\zeta(s-1)}{\zeta(s)} \cdot \frac{x^s}{s^2} ds.
\end{equation}
The integrand above has a simple pole at $s=2$ with residue
$\frac{1}{\zeta(2)}\cdot \frac{x^2}{4}=\frac{3x^2}{2\pi^2}$. Thus,
by moving the integration contour to $\Real s=1$ and using the
trivial estimate
\begin{equation*}
\int\limits_1^{\sigma_0} \frac{\vert x^{s\pm iR}\vert}{\vert
s+iR\vert^2} ds \ll \frac{x^{\sigma_0}}{R},
\end{equation*}
we infer from \eqref{4.15} the equality
\begin{equation}\label{4.16}
\sum\limits_{1\leq q<x} \varphi (q) \log \frac{x}{q}
=\frac{3x^2}{\pi^2} +\frac{1}{2\pi i}
\int\limits_{1-i\infty}^{1+i\infty} \frac{\zeta(s-1)}{\zeta(s)}
\cdot \frac{x^s}{s^2} ds.
\end{equation}
Utilizing the functional equation
\begin{equation*}
\frac{\zeta(it)}{\zeta(1+it)} =\chi (it)\cdot
\frac{\zeta(1-it)}{\zeta(1+it)}, \qquad \chi
(z)=\frac{(2\pi)^z}{2\Gamma (z)\cos \frac{\pi z}{2}},
\end{equation*}
and the equalities $\zeta (1-it)=\overline{\zeta (1+it)}$ and (see
relation 8.332 in \cite{GR})
\begin{equation*}
\vert \Gamma (it)\vert^2 =\frac{\pi}{t\sinh \pi t}\qquad (t>0),
\end{equation*}
we find that the expression in \eqref{4.16} is given by
\begin{equation*}
\frac{3x^2}{2\pi^2}+O\left( \int\limits_0^\infty \frac{x\sqrt{t}}{
(1+t^2)} dt\right) =\frac{3x^2}{2\pi^2}+O(x),
\end{equation*}
which completes the proof of \eqref{1.6}.

Note that one gets a better error term than in \eqref{2.2} as a
result of the presence of the factor $\log \frac{x}{q}$ instead of
$1$.

\bigskip

{\bf Acknowledgements.} We are grateful to Richard R. Hall for the
argument leading to the proof of \eqref{1.6}, to Kevin Ford for
bringing to our attention reference \cite{Wal}, and to the referee
for several useful suggestions.


\bigskip

\bigskip

\begin{thebibliography}{99}

\bibitem{ABCZ}
V. Augustin, F. P. Boca, C. Cobeli, A. Zaharescu, `The $h$-spacing
distribution between Farey points', {\em Math.\ Proc.\ Camb.\
Phil.\ Soc.\ }131 (2001), 23--38.
%
\bibitem{BCZ0}
F. P. Boca, C. Cobeli, A. Zaharescu, `Distribution of lattice
points visible from the origin', {\em Comm.\ Math.\ Phys.\ }213
(2000), 433--470.
%
\bibitem{BCZ}
F. P. Boca, C. Cobeli, A. Zaharescu, `A conjecture of R.R. Hall on
Farey points', {\em J.\ Reine\ Angew.\ Mathematik\ }535 (2001),
207--236.
%
\bibitem{BGZ}
F. P. Boca, R. N. Gologan, A. Zaharescu, `The average length of a
trajectory in a certain billiard in a flat two-torus', {\em New\
York\ J.\ Math.\ }9 (2003), 303--330.
%
\bibitem{BZ}
F. P. Boca, A. Zaharescu, `The distribution of the free path
lengths in the periodic two-dimensional Lorentz gas in the
small-scatterer limit', preprint math.NT/0301270.
%
\bibitem{Edw}
H. M. Edwards, {\em Riemann's zeta function}, (Dover Publications,
Inc., Mineola, New York, 2001).
%
\bibitem{Fra}
J. Franel, `Les suites de Farey et le probl\` eme des nombres
premiers', {\em Gottinger\ Nachr.\ }(1924), 198--201.
%
\bibitem{GR}
I. S. Gradshteyn, I. M. Ryzhik, {\em Table of integrals, series,
and products}, (6th edition, Academic Press, 2000).
%
\bibitem{Hall}
R. R. Hall, `A note on Farey series', {\em J.\ London\ Math.\
Soc.\ }2 (1970), 139--148.
%
\bibitem{Hej}
D. A. Hejhal, `On the triple correlations of zeros of the zeta
function', {\em Internat.\ Math.\ Res.\ Notices\ }1994, no. 7.
%
\bibitem{HZ}
M. N. Huxley, A. Zhigljavsky, `On the distribution of Farey
fractions and hyperbolic lattice points', {\em Period.\ Math.\
Hungarica\ }42 (2001), 191--198.
%
\bibitem{KS}
N. M. Katz, P. Sarnak, `Zeros of zeta functions and symmetry',
{\em Bull.\ Amer.\ Math.\ Soc.\ }36 (1999), 1--26.
%
\bibitem{Lan}
E. Landau,  `Bemerkungen zu der vorstehenden Abhandlung von Herrn
Franel', {\em Gottinger\ Nachr.\ }(1924), 202--206.
%
\bibitem{Mo}
H. L. Montgomery, `The pair correlation of zeros of the zeta
function', in {\em Analytic number theory} (Proc. Sympos. Pure
Math., Vol. XXIV, St. Louis Univ., Mo., 1972), pp. 181--193, Amer.
Math. Soc., Providence, R.I., 1973.
%
\bibitem{RS}
Z. Rudnick, P. Sarnak, `Zeros of principal L-functions and random
matrix theory', {\em Duke\ Math.\ J.\ } 81 (1996), 269--322.
%
\bibitem{Sar}
P. Sarnak, {\em Some applications of modular forms}, (Cambridge
Tracts in Mathematics 99, Cambridge University Press, 1990).
%
\bibitem{Wal}
A. Walfisz, {\em Weylsche Exponentialsummen in der neueren
Zahlentheorie} (VEB Deutscher Verlag der Wissenschaften, Berlin,
1963).
\end{thebibliography}
\end{document}